\documentclass{amsart}
\usepackage{amssymb}

\usepackage{graphicx}

\input{tcilatex}
\input tcilatex

\begin{document}
\title[Duality and Discrete Non-Neutral Wright-Fisher Models]{A Duality Approach to the Genealogies of Discrete Non-Neutral Wright-Fisher
Models}
\author{Thierry E. Huillet}
\address{Laboratoire de Physique Th\'{e}orique et Mod\'{e}lisation \\
CNRS-UMR 8089 et Universit\'{e} de Cergy-Pontoise, 2 Avenue Adolphe Chauvin,
95302, Cergy-Pontoise, FRANCE\\
E-mail: Thierry.Huillet@u-cergy.fr}
\maketitle

\begin{abstract}
Discrete ancestral problems arising in population genetics are investigated.
In the neutral case, the duality concept has proved of particular interest
in the understanding of backward in time ancestral process from the forward
in time branching population dynamics. We show that duality formulae still
are of great use when considering discrete non-neutral Wright-Fisher models.
This concerns a large class of non-neutral models with completely monotone
(CM) bias probabilities. We show that most classical bias probabilities used
in the genetics literature fall within this CM class or are amenable to it
through some `reciprocal mechanism' which we define. Next, using elementary
algebra on CM functions, some suggested novel evolutionary mechanisms of
potential interest are introduced and discussed, \newline

\textbf{Running title:} Duality and Discrete Non-Neutral Wright-Fisher
Models.\newline

\end{abstract}

\section{Introduction}

Forward evolution of large populations in genetics has a long history,
starting in the $1920$s; it is closely attached to the names of R. A. Fisher
and S. Wright; see T. Nagylaki ('$1999$) for historical commentaries and on
the role played by the French geneticist G. Mal\'{e}cot, starting shortly
before the second world war. The book of W. Ewens ('$2004$) is an excellent
modern presentation of the current mathematical theory. Coalescent theory is
the corresponding backward problem, obtained while running the forward
evolution processes backward-in-time. It was discovered independently by
several researchers in the $1980$s, but definitive formalization is commonly
attributed to J. Kingman ('$1982$). Major contributions to the development
of coalescent theory were made (among others) by P. Donnelly, R. Griffiths,
R. Hudson, F. Tajima and S. Tavar\'{e} (see the course of Tavar\'{e} in
Saint-Flour '$2004$ for a review). It included incorporating variations in
population size, mutation, recombination, selection... In ('$1999$), J.
Pitman and S. Sagitov, independently, introduced coalescent processes with
multiple collisions of ancestral lineages. Shortly later, the full class of
exchangeable coalescent processes with simultaneous multiple mergers of
ancestral lineages was discovered by M. M\"{o}hle and S. Sagitov ('$2001$)
and J. Schweinsberg ('$2000$). All these recent developments and
improvements concern chiefly the discrete neutral case and their various
scaling limits in continuous time and/or space. As was shown by M\"{o}hle ('$%
1994$ and '$1999$), neutral forward and backward theories learn much from
one another by using a concept of duality introduced by T. Liggett ('$1985$%
). Backward theory in the presence of mutations in the forward process is
well-understood, as it requires the study of a marked Kingman's tree (see
Tavar\'{e}, ('$2004$) for a review). In the works of C. Neuhauser and S.
Krone ('$1997$), there is also some use of the duality concept in an attempt
to understand the genealogies of a Wright-Fisher diffusion (as a limit of a
discrete Wright-Fisher model) presenting a selection mechanism; this led
these authors to the idea of the ancestral selection graph extending
Kingman's coalescent tree of the neutral theory; see also T. Huillet ('$2007$%
) for related objectives in the context of Wright-Fisher diffusions with and
without drifts. There is therefore some evidence that the concept of duality
could help one understand the backward theory even in non-neutral situations
when various evolutionary forces are the causes of deviation to neutrality
(see J. Crow and M. Kimura, ('$1970$), T. Maruyama ('$1977$), J. Gillepsie ('%
$1991$) and W. Ewens ('$2004$), for a discussion on various models of utmost
interest in population genetics).

In this Note, we focus on discrete non-neutral Wright-Fisher (say WF) models
and on the conditions on the bias probabilities under which forward
branching population dynamics is directly amenable to a dual discrete
ancestral coalescent. We emphasize that duality formulae still are of great
use when considering discrete non-neutral Wright-Fisher models, at least for
specific deviation forces to neutrality. It is shown that it concerns a
large class of non-neutral models involving completely monotone bias
probabilities. Several classical examples are supplied in the light of
complete monotonicity. In the process leading us to focus on these peculiar
bias models, some unsuspected evolutionary mechanisms of potential interest
are introduced and discussed, as suggested by elementary algebra on
completely monotone functions. We emphasize that the relevance of these
novel bias mechanisms in Biology seems to deserve additional work and
confrontation with real-world problems is urged for to pinpoint their
biological significance.

We shall finally briefly outline the content of this manuscript. Section $2$
is designed to fix the background and ideas: We introduce some basic facts
about the discrete-time forward (subsection $2.2$) and backward processes
(subsection $2.3$) arising from exchangeable reproduction laws (subsection $%
2.1$). In subsection $2.4$, we introduce a concept of duality and briefly
recall its relevance to the study of the neutral case problem. The basic
question we address in subsequent sections is whether this notion of duality
still makes sense in non-neutral situations. We start supplying important
non-neutral examples in section $3$. In section $4$, we show that duality
does indeed make sense in the framework of discrete non-neutral
Wright-Fisher models, but only for the class of completely-monotone
state-dependent transition frequencies. In section $5$, we show that most
non-neutrality mechanisms used in the literature fall within this class, or
are amenable to it via some `reciprocal transformation', starting with
elementary mechanisms and ending up with more complex ones. In section $6,$
we show that duality can be used in non-neutral situations to compute the
extinction probabilities (invariant measure) of the dual backward ancestral
process\emph{\ }if one knows the invariant measure (respectively, extinction
probabilities) of the forward branching process.

\section{Discrete-time neutral coalescent}

In this Section, to fix the background and notations, we review some
well-known facts from the cited literature.

\subsection{Exchangeable neutral population models: Reproduction laws.}

(The Cannings model: '$1974$). Consider a population with non-overlapping
generations $r\in \Bbb{Z}.$ Assume the population size is constant, say $n$ (%
$n$ individuals (or genes)) over generations. Assume the random reproduction
law at generation $0$ is $\mathbf{\nu }_{n}:=\left( \nu _{1,n},...,\nu
_{n,n}\right) ,$ satisfying: 
\begin{equation*}
\sum_{m=1}^{n}\nu _{m,n}=n.
\end{equation*}
Here, $\nu _{m,n}$ is the number of offspring of gene $m.$ We avoid the
trivial case: $\nu _{m,n}=1$, $m=1,...,n.$ One iterates the reproduction
over generations, while imposing the following additional assumptions:

- Exchangeability: $\left( \nu _{1,n},...,\nu _{n,n}\right) \overset{d}{=}%
\left( \nu _{\sigma \left( 1\right) ,n},...,\nu _{\sigma \left( n\right)
,n}\right) ,$ for all permutations $\sigma \in \mathcal{S}_{n}.$

- time-homogeneity: reproduction laws are independent and identically
distributed (iid) at each generation $r\in \Bbb{Z}.$

This model therefore consists of a conservative conditioned branching
Galton-Watson process in $\left[ n\right] ^{\Bbb{Z}}$, where $\left[
n\right] :=\left\{ 0,1,...,n\right\} $ (see Karlin-McGregor, '$1964$).%
\newline

Famous reproduction laws are:

\emph{Example 2.1.1} The multinomial Dirichlet model: $\mathbf{\nu }_{n}%
\overset{d}{\sim }$ Multin-Dirichlet$\left( n;\theta \right) $, where $%
\theta >0$ is a disorder parameter. With $\mathbf{k}_{n}:=\left(
k_{1},...,k_{n}\right) $, $\mathbf{\nu }_{n}$ admits the following joint
exchangeable distribution on the simplex $\left| \mathbf{k}_{n}\right|
:=\sum_{m=1}^{n}k_{m}=n$: 
\begin{equation*}
\Bbb{P}_{\theta }\left( \mathbf{\nu }_{n}=\mathbf{k}_{n}\right) =\frac{n!}{%
\left[ n\theta \right] _{n}}\prod_{m=1}^{n}\frac{\left[ \theta \right]
_{k_{m}}}{k_{m}!},
\end{equation*}
where $\left[ \theta \right] _{k}=\theta \left( \theta +1\right) ...\left(
\theta +k-1\right) $ is the rising factorial of $\theta $. This distribution
can be obtained by conditioning $n$ independent mean $1$ P\`{o}lya
distributed random variables $\mathbf{\xi }_{n}=\left( \xi _{1},...,\xi
_{n}\right) $ on summing to $n$, that is to say: $\mathbf{\nu }_{n}\overset{d%
}{=}\left( \mathbf{\xi }_{n}\mid \left| \mathbf{\xi }_{n}\right| =n\right) ,$
where 
\begin{equation*}
\Bbb{P}_{\theta }\left( \xi _{1}=k\right) =\frac{\left[ \theta \right] _{k}}{%
k!}\left( 1+\theta \right) ^{-k}\left( \theta /\left( 1+\theta \right)
\right) ^{\theta }\text{, }k\in \Bbb{N}.
\end{equation*}

When $\theta \uparrow \infty $, this distribution reduces to the
Wright-Fisher model for which $\mathbf{\nu }_{n}\overset{d}{\sim }$ Multin$%
\left( n;1/n,...,1/n\right) .$ Indeed, $\mathbf{\nu }_{n}$ admits the
following joint exchangeable multinomial distribution on the simplex $\left| 
\mathbf{k}_{n}\right| =n$: 
\begin{equation*}
\Bbb{P}_{\infty }\left( \mathbf{\nu }_{n}=\mathbf{k}_{n}\right) =\frac{%
n!\cdot n^{-n}}{\prod_{m=1}^{n}k_{m}!}.
\end{equation*}
This distribution can be obtained by conditioning $n$ independent mean $1$
Poisson distributed random variables $\mathbf{\xi }_{n}=\left( \xi
_{1},...,\xi _{n}\right) $ on summing to $n$: $\mathbf{\nu }_{n}\overset{d}{=%
}\left( \mathbf{\xi }_{n}\mid \left| \mathbf{\xi }_{n}\right| =n\right) $.
When $n$ is large, using Stirling formula, $n!\sim \sqrt{2\pi }%
n^{n+1/2}e^{-n};$ it follows that $\mathbf{\nu }_{n}\overset{d}{\underset{%
n\uparrow \infty }{\rightarrow }}\mathbf{\xi }_{\infty }$ with joint
finite-dimensional law: $\Bbb{P}\left( \mathbf{\xi }_{n}=\mathbf{k}%
_{n}\right) =\prod_{m=1}^{n}\frac{e^{-1}}{k_{m}!}=\frac{e^{-n}}{%
\prod_{m=1}^{n}k_{m}!}$ on $\Bbb{N}^{n}.$ Thanks to the product form of all
finite-dimensional laws of $\mathbf{\xi }_{\infty }$, we get an asymptotic
independence property of $\mathbf{\nu }_{n}$.

\emph{Example 2.1.2} In the Moran model, $\mathbf{\nu }_{n}\overset{d}{\sim }
$ random permutation of $\left( 2,0,1,...,1\right) :$ in such a model, only
one new gene per generation may come to life, at the expense of the
simultaneous disappearance of some other gene.

\subsection{Forward in time branching process}

Take a sub-sample of size $m$ from $\left[ n\right] :=\left\{
0,1,...,n\right\} $ at generation $0.$ Let

\begin{equation*}
N_{r}\left( m\right) \text{ }=\#\text{ offspring at generation }r\in \Bbb{N}%
_{+}\text{, forward-in-time}.
\end{equation*}
This sibship process is a discrete-time homogeneous Markov chain, with
transition probability: 
\begin{equation}
\Bbb{P}\left( N_{r+1}\left( m\right) =k^{\prime }\mid N_{r}\left( m\right)
=k\right) =\Bbb{P}\left( \nu _{1,n}+...+\nu _{k,n}=k^{\prime }\right) .
\label{Eq1}
\end{equation}
It is a martingale, with state-space $\left\{ 0,...,n\right\} $, initial
state $m$, absorbing states $\left\{ 0,n\right\} $ and transient states $%
\left\{ 1,...,n-1\right\} .$ The first hitting time of boundaries $\left\{
0,n\right\} $, which is: $\tau \left( m\right) =\tau _{\left\{ 0\right\}
}\left( m\right) \wedge \tau _{\left\{ n\right\} }\left( m\right) $ is
finite with probability $1$ and has finite mean. Omitting reference to any
specific initial condition $m$, the process $\left( N_{r};r\in \Bbb{N}%
\right) $ has the transition matrix $\Pi _{n}$ with entries $\Pi _{n}\left(
k,k^{\prime }\right) =\Bbb{P}\left( \nu _{1,n}+...+\nu _{k,n}=k^{\prime
}\right) $ given by (\ref{Eq1})$.$ We have $\Pi _{n}\left( 0,k^{\prime
}\right) =\delta _{0,k^{\prime }}$ and $\Pi _{n}\left( n,k^{\prime }\right)
=\delta _{n,k^{\prime }}$ and $\Pi _{n}$ is not irreducible. However, $\Pi
_{n}$ is aperiodic and (apart from absorbing states) cannot be broken down
into non-communicating subsets; as a result it is diagonalizable, with
eigenvalues $\left| \lambda _{0}\right| \geq \left| \lambda _{1}\right| \geq
...\geq \left| \lambda _{n}\right| $ and $1=\lambda _{0}=\lambda _{1}>\left|
\lambda _{2}\right| $.\newline

\emph{Example 2.2.1} (Dirichlet binomial): With $U$ a $\left( 0,1\right) -$%
valued random variable with density beta$\left( k\theta ,\left( n-k\right)
\theta \right) $%
\begin{equation*}
\Bbb{P}\left( \nu _{1,n}+..+\nu _{k,n}=k^{\prime }\right) =\binom{n}{%
k^{\prime }}\frac{\left[ k\theta \right] _{k^{\prime }}\left[ \left(
n-k\right) \theta \right] _{n-k^{\prime }}}{\left[ n\theta \right] _{n}}=%
\Bbb{E}\left[ \binom{n}{k^{\prime }}U^{k^{\prime }}\left( 1-U\right)
^{n-k^{\prime }}\right] ,
\end{equation*}
which is a beta mixture of the binomial distribution Bin$\left( n,u\right) .$

\emph{Example 2.2.2} The Wright-Fisher model has a Bin$\left( n,k/n\right) $
transition matrix: 
\begin{equation*}
\Bbb{P}\left( N_{r+1}\left( m\right) =k^{\prime }\mid N_{r}\left( m\right)
=k\right) =\binom{n}{k^{\prime }}\left( \frac{k}{n}\right) ^{k^{\prime
}}\left( 1-\frac{k}{n}\right) ^{n-k^{\prime }}.
\end{equation*}
\textbf{Remark (}statistical symmetry\textbf{)}: Due to exchangeability of
the reproduction law, neutral models are symmetric in the following sense:
The transition probabilities of $\overline{N}_{r}\left( m\right)
:=n-N_{r}\left( m\right) $ are equal to the transition probabilities of $%
N_{r}\left( m\right) $. $\square $

\subsection{Backward in time process}

\textbf{(neutral coalescent)}

$\emph{The}$ \emph{c}$\emph{oalescent}$ $\emph{backward}$ $\emph{process}$
can be defined as follows: Take a sub-sample of size $m$ from $\left[
n\right] $ at generation $0.$ Identify two individuals from $\left[ m\right] 
$ at each step if they share a common ancestor one generation
backward-in-time. This defines an equivalence relation between $2$ genes
from the set $\left[ m\right] $. Define the induced ancestral backward
process:

\begin{center}
$\mathcal{A}_{r}\left( m\right) $ $\in \mathcal{E}_{m}=\left\{ \text{%
equivalence classes (partitions) of }\left[ m\right] \right\} ,$ $r\in \Bbb{N%
}$, backward-in-time$.$
\end{center}

The ancestral process is a discrete-time Markov chain with transition
probability: 
\begin{equation*}
\Bbb{P}\left( \mathcal{A}_{r+1}\left( m\right) =\alpha \mid \mathcal{A}%
_{r}\left( m\right) =\beta \right) =P_{\beta ;\alpha }\text{; with }\left(
\alpha ,\beta \right) \in \mathcal{E}_{m}\text{, }\alpha \subseteq \beta
\end{equation*}
where, with $a=\left| \alpha \right| =$ number of equivalence classes of $%
\alpha ,$ $b=\left| \beta \right| =$ number of equivalence classes of $\beta
,$ $\mathbf{b}_{a}:=\left( b_{1},...,b_{a}\right) $ clusters sizes of $\beta 
$ and $\left( m\right) _{a}:=m\left( m-1\right) ...\left( m-a+1\right) $ a
falling factorial,

\begin{equation*}
P_{\beta ;\alpha }=P_{b;a}^{\left( n\right) }\left( \mathbf{b}_{a}\right) =%
\frac{\left( n\right) _{a}}{\left( n\right) _{b}}\Bbb{E}\left(
\prod_{l=1}^{a}\left( \nu _{l,n}\right) _{b_{l}}\right) \text{ }
\end{equation*}
is the probability of a $\mathbf{b}_{a}-$merger. This is the probability
that $b$ randomly chosen individuals out of $n$ have $a\leq b$ distinct
parents, $c$ merging classes and cluster sizes $b_{1}\geq ...\geq b_{c}\geq
2 $, $b_{c+1}=...=b_{a}=1.$

If $c=1$: a unique multiple collision occurs of order $b_{1}\geq 2$.

If $b_{1}=2$: a simple binary collision occurs involving only two clusters.

If $c>1$, simultaneous multiple collisions of orders $b_{1}\geq ...\geq
b_{c}\geq 2$ occur.

Thus, the jump's height of a transition $b\rightarrow a$ is $%
b-a=\sum_{i=1}^{c}\left( b_{i}-1\right) ,$ corresponding to a partition of
integer $b-a$ into $c$ summands, each $\geq 1$.

The chain's state-space is: $\left\{ \text{equivalence relations on
(partitions of) }\left\{ 1,...,m\right\} \right\} $; it has dimension $%
B_{m}:=\sum_{l=0}^{m}$ $S_{m,l}$ (a Bell number), where $S_{m,l}$ are the
second-kind Stirling numbers.

The chain has initial state $\mathcal{A}_{0}=\left\{ \left( 1\right)
,...,\left( m\right) \right\} $, and terminal absorbing state $\left\{
\left( 1,...,m\right) \right\} .$

\emph{Examples: }

From the Dirichlet example 2.2.1, we get: $P_{b;a}^{\left( n\right) }\left( 
\mathbf{b}_{a}\right) =\frac{\left( n\right) _{a}}{\left[ n\theta \right]
_{b}}\prod_{m=1}^{a}\left[ \theta \right] _{b_{m}}$.

From the WF example 2.2.2: In this case, $P_{b;a}^{\left( n\right) }\left( 
\mathbf{b}_{a}\right) =\frac{\left( n\right) _{a}}{n^{b}}$ is the uniform
distribution on $\left\{ \mathbf{b}_{a}:b_{1}+...+b_{a}=b\right\} $. $%
\newline
$

\emph{The ancestral Count Process:} Let

\begin{equation*}
A_{r}\left( m\right) \text{ }=\#\text{ ancestors at generation }r\in \Bbb{N}%
\text{, backward-in-time}.
\end{equation*}
\begin{equation*}
\text{Then: }A_{r}\left( m\right) =\#\text{ blocks of }\mathcal{A}_{r}\left(
m\right) .
\end{equation*}
The backward ancestral count process is a discrete-time Markov chain with
transition probabilities (Cannings, '$1974$ and Gladstien '$1978$): 
\begin{equation}
\Bbb{P}\left( A_{r+1}\left( m\right) =a\mid A_{r}\left( m\right) =b\right)
=P_{b,a}^{\left( n\right) }:=\frac{b!}{a!}\sum_{\underset{b_{1}+...+b_{a}=b}{%
b_{1},...,b_{a}\in \Bbb{N}_{+}}}\frac{P_{b;a}^{\left( n\right) }\left( 
\mathbf{b}_{a}\right) }{b_{1}!...b_{a}!}.  \label{Eq2}
\end{equation}
\begin{equation*}
=\frac{\binom{n}{a}}{\binom{n}{b}}\sum_{\underset{b_{1}+...+b_{a}=b}{%
b_{1},...,b_{a}\in \Bbb{N}_{+}}}\Bbb{E}\left( \prod_{l=1}^{a}\binom{\nu
_{l,n}}{b_{l}}\right) .
\end{equation*}
This Markov chain has state-space $\left\{ 0,...,m\right\} $, initial state $%
m$, absorbing states $\left\{ 0,1\right\} .$ The process $\left( A_{r};r\in 
\Bbb{N}\right) $ has the transition matrix $P_{n}$ with entries $P_{n}\left(
b,a\right) =P_{b,a}^{\left( n\right) }$ given by (\ref{Eq2})$.$ Note, by
inclusion-exclusion principle, the alternative alternating expression of $%
P_{b,a}^{\left( n\right) }:$%
\begin{equation*}
P_{b,a}^{\left( n\right) }:=\frac{\binom{n}{a}}{\binom{n}{b}}%
\sum_{m=0}^{a}\left( -1\right) ^{a-m}\binom{a}{m}\Bbb{E}\left( \binom{\nu
_{1,n}+...+\nu _{m,n}}{b}\right) .
\end{equation*}

\subsection{Duality (neutral case).}

We start with a definition of the duality concept which is relevant to our
purposes.

\textbf{Definition }[Liggett, '$1985$]: \emph{Two Markov processes }$\left(
X_{t}^{1},X_{t}^{2};t\geq 0\right) ,$\emph{\ with state-spaces }$\left( 
\mathcal{E}_{1},\mathcal{E}_{2}\right) ,$\emph{\ are said to be dual with
respect to some real-valued function }$\Phi $\emph{\ on the product space }$%
\mathcal{E}_{1}\times \mathcal{E}_{2}$\emph{\ if }$\forall x_{1}\in \mathcal{%
E}_{1},$\emph{\ }$\forall x_{2}\in \mathcal{E}_{2},$\emph{\ }$\forall t\geq
0:$%
\begin{equation}
\Bbb{E}_{x_{1}}\Phi \left( X_{t}^{1},x_{2}\right) =\Bbb{E}_{x_{2}}\Phi
\left( x_{1},X_{t}^{2}\right) .  \label{dual}
\end{equation}
\newline
We then recall basic examples of dual processes from the neutral and
exchangeable population models (M\"{o}hle, '$1997$): The neutral forward and
backward processes $\left( N_{r},A_{r};r\in \Bbb{N}\right) $ introduced in
the two preceding subsections are dual with respect to the hypergeometric
sampling without replacement kernels: 
\begin{equation}
\left( i\text{ }\right) \text{ }\Phi _{n}^{1}\left( m,k\right) =\binom{m}{k}/%
\binom{n}{k}\text{ and }  \label{Eq3}
\end{equation}
\begin{equation*}
\left( ii\right) \text{ }\Phi _{n}^{2}\left( m,k\right) =\binom{n-m}{k}/%
\binom{n}{k}\text{ on }\left\{ 0,...,n\right\} ^{2}.
\end{equation*}
Namely $\left( i\right) $ reads: 
\begin{equation*}
\text{ }\Bbb{E}_{m}\left[ \binom{N_{r}}{k}/\binom{n}{k}\right] =\Bbb{E}%
_{k}\left[ \binom{m}{A_{r}}/\binom{n}{A_{r}}\right] =\Bbb{E}_{k}\left[ 
\binom{n-A_{r}}{n-m}/\binom{n}{n-m}\right] .
\end{equation*}
Call type $A$ individuals the descendants of the $m$ first chosen
individuals (allele $A$) in the study of the forward process; type $a$
individuals are the remaining ones (allele $a$). The left-hand-side of the
above equality is an expression of the probability that a $k-$sample
(without replacement) from population of size $N_{r}$ at time $r$ are all of
type $A$, given $N_{0}=m.$ If this $k-$sample are all descendants of $A_{r}$
ancestors at time $-r$, this probability must be equal to the probability
that a $\left( n-m\right) -$sample from population of size $A_{r}$ at time $%
-r$ are all of type $a$. This is the meaning of the right-hand-side.

And $\left( ii\right) $ reads: 
\begin{equation*}
\text{ }\Bbb{E}_{m}\left[ \binom{n-N_{r}}{k}/\binom{n}{k}\right] =\Bbb{E}%
_{k}\left[ \binom{n-m}{A_{r}}/\binom{n}{A_{r}}\right] =\Bbb{E}_{k}\left[ 
\binom{n-A_{r}}{m}/\binom{n}{m}\right] .
\end{equation*}
The left-hand-side is the probability that a $k-$sample (without
replacement) from population of size $N_{r}$ at time $r$ are all of type $a$%
, given $N_{0}=m.$ If this $k-$sample are all descendants of $A_{r}$
ancestors at time $-r$, this probability must be equal to the probability
that a $m-$sample from population of size $A_{r}$ at time $-r$ are
themselves all of type $a.$\newline

With $P_{n}^{\prime }$ the transpose of $P_{n}$, a one-step ($r=1$) version
of these formulae is: 
\begin{equation*}
\left( i\right) \text{ }\Pi _{n}\Phi _{n}^{1}=\Phi _{n}^{1}P_{n}^{\prime }%
\text{ and }\left( ii\right) \text{ }\Pi _{n}\Phi _{n}^{2}=\Phi
_{n}^{2}P_{n}^{\prime }
\end{equation*}
where $\left( \Phi _{n}^{1},\Phi _{n}^{2}\right) $ are $n\times n$ matrices
with entries $\Phi _{n}^{1}\left( m,k\right) $ and $\Phi _{n}^{2}\left(
m,k\right) ,$ respectively and $\left( \Pi _{n},P_{n}\right) $ the
transition matrices of forward and backward processes. Note that the matrix $%
\Phi _{n}^{2}$ is symmetric and left-upper triangular. The matrices $\Phi
_{n}^{1}$ and $\Phi _{n}^{2}$ are both invertible, with respective entries 
\begin{equation*}
\text{ }\left[ \Phi _{n}^{1}\right] ^{-1}\left( i,j\right) =\left( -1\right)
^{i-j}\binom{i}{j}\binom{n}{i}
\end{equation*}
and 
\begin{equation*}
\text{ }\left[ \Phi _{n}^{2}\right] ^{-1}\left( i,j\right) =\left( -1\right)
^{i+j-n}\binom{i}{n-j}\binom{n}{i}=\left( -1\right) ^{i+j-n}\binom{j}{n-i}%
\binom{n}{j}.
\end{equation*}
The matrix $\left[ \Phi _{n}^{1}\right] ^{-1}$ is left-lower triangular,
while $\left[ \Phi _{n}^{2}\right] ^{-1}$ is symmetric right-lower
triangular. Thus, 
\begin{equation*}
\left( i\right) \text{ }\left[ \Phi _{n}^{1}\right] ^{-1}\Pi _{n}\Phi
_{n}^{1}=P_{n}^{\prime }\text{ and }\left( ii\right) \text{ }\left[ \Phi
_{n}^{2}\right] ^{-1}\Pi _{n}\Phi _{n}^{2}=P_{n}^{\prime }.
\end{equation*}
In any case, being similar matrices, $\Pi _{n}$ and $P_{n}^{\prime }$ (or $%
P_{n}$) both share the same eigenvalues. If $R_{n}$ diagonalizing $\Pi _{n}$
is known so that $R_{n}^{-1}\Pi _{n}R_{n}=\Lambda _{n}:=$ diag$\left(
\lambda _{0},...,\lambda _{n}\right) ,$ the diagonal matrix of the
eigenvalues of $\Pi _{n}$, then, with $\Phi _{n}=\Phi _{n}^{1}$ or $\Phi
_{n}^{2}$, $\widetilde{R}_{n}:=\Phi _{n}^{-1}R_{n}$ diagonalizes $%
P_{n}^{\prime }$ and is obtained for free (and conversely). $R_{n}$ is the
matrix whose columns are the right-eigenvectors of $\Pi _{n}$ and $%
\widetilde{R}_{n}$ is the matrix whose columns (rows) are the
right-eigenvectors (left-eigenvectors) of $P_{n}^{\prime }$ (of $P_{n}$).
Similarly, if $L_{n}$ is the matrix whose rows are the left-eigenvectors of $%
\Pi _{n}$, $\widetilde{L}_{n}:=L_{n}\Phi _{n}$ is the matrix whose rows
(columns) are the left-eigenvectors (right-eigenvectors) of $P_{n}^{\prime }$
(of $P_{n}$). With $l_{k}^{\prime }$ the $k-$th row of $L_{n}$ and $r_{k}$
the $k-$th column of $R_{n},$ the spectral decomposition of $\Pi _{n}$ is: 
\begin{equation*}
\Pi _{n}^{r}=\sum_{k=0}^{n}\lambda _{k}^{r}\frac{r_{k}l_{k}^{\prime }}{%
l_{k}^{\prime }r_{k}}\text{, }r\in \Bbb{N},
\end{equation*}
whereas, with $\widetilde{l}_{k}$ the $k-$th column of $\widetilde{L}_{n}$
and $\widetilde{r}_{k}^{\prime }$ the $k-$th row of $\widetilde{R}_{n},$ the
one of $P_{n}$ reads: 
\begin{equation*}
P_{n}^{r}=\sum_{k=0}^{n}\lambda _{k}^{r}\frac{\widetilde{l}_{k}\widetilde{r}%
_{k}^{\prime }}{\widetilde{r}_{k}^{\prime }\widetilde{l}_{k}}%
=\sum_{k=0}^{n}\lambda _{k}^{r}\frac{\Phi _{n}^{\prime }l_{k}\left( \Phi
_{n}^{-1}r_{k}\right) ^{\prime }}{\left( \Phi _{n}^{-1}r_{k}\right) ^{\prime
}\Phi _{n}^{\prime }l_{k}}\text{, }r\in \Bbb{N}.
\end{equation*}
\newline

In M\"{o}hle '$1999$, a direct combinatorial proof of the duality result can
be found (in the general exchangeable or neutral case); it was obtained by
directly checking the consistency of (\ref{Eq1}), (\ref{Eq2}) and (\ref{Eq3}%
).

The duality formulae allow one to deduce the probabilistic structure of one
process from the one of the other. The question we address now is: does
duality still make sense in non-neutral situations? We shall see that it
does in discrete non-neutral Wright-Fisher models, but only for some class
of state-dependent transition frequencies.

\section{Beyond neutrality (symmetry breaking)}

Discrete forward non-neutral models (with non-null drifts) can be obtained
by substituting 
\begin{equation*}
k\rightarrow np\left( \frac{k}{n}\right) \text{ in }\Bbb{P}\left( \nu
_{1,n}+...+\nu _{k,n}=k^{\prime }\right) ,
\end{equation*}
where: 
\begin{equation*}
p\left( x\right) :x\in \left( 0,1\right) \rightarrow \left( 0,1\right) \text{
is continuous, increasing, with }p\left( 0\right) =0,\text{ }p\left(
1\right) =1.
\end{equation*}
$p\left( x\right) $ is the state-dependent Bernoulli bias probability
different from identity $x$ (as in neutral case). \newline

When particularized to the WF model, this leads to the biased transition
probabilities: 
\begin{equation*}
\Bbb{P}\left( N_{r+1}\left( m\right) =k^{\prime }\mid N_{r}\left( m\right)
=k\right) =\binom{n}{k^{\prime }}p\left( \frac{k}{n}\right) ^{k^{\prime
}}\left( 1-p\left( \frac{k}{n}\right) \right) ^{n-k^{\prime }}.
\end{equation*}
In this binomial $n-$sampling with replacement model, a type $A$ individual
is drawn with probability $p\left( \frac{k}{n}\right) $ which is different
from the uniform distribution $k/n$, due to bias effects$.$

From this, we conclude (a symmetry breaking property): The transition
probabilities of $\overline{N}_{r}\left( m\right) :=n-N_{r}\left( m\right) $%
, $r\in \Bbb{N}$ are 
\begin{equation*}
\text{Bin}\left( n,1-p\left( 1-k/n\right) \right) \neq \text{Bin}\left(
n,p\left( k/n\right) \right) ,
\end{equation*}
and so, no longer coincide with the ones of $\left( N_{r}\left( m\right)
;r\in \Bbb{N}\right) .$ The process $N_{r}\left( m\right) $, $r\in \Bbb{N}$
no longer is a martingale. Rather, if $x\rightarrow p\left( x\right) $ is
concave (convex), $N_{r}\left( m\right) $, $r\in \Bbb{N}$ is a submartingale
(supermartingale), because: $\Bbb{E}\left( N_{r+1}\left( m\right) \mid
N_{r}\left( m\right) \right) =np\left( N_{r}\left( m\right) /n\right) \geq
N_{r}\left( m\right) $ (respectively $\leq N_{r}\left( m\right) $)$.$

In the binomial neutral Wright-Fisher transition probabilities, we replaced
the success probability $\frac{k}{n}$ by a more general function $p\left( 
\frac{k}{n}\right) $. However, this replacement leaves open the question
what model is in the background and what quantity the process $\left(
N_{r},r\in \Bbb{N}\right) $ really counts. A concrete model in terms of
offspring variables must be provided instead. To address this question, we
emphasize that the reproduction law corresponding to the biased binomial
model is multinomial and asymmetric, namely: $\mathbf{\nu }_{n}\overset{d}{%
\sim }$ Multin$\left( n;\mathbf{\pi }_{n}\right) $, where $\mathbf{\pi }%
_{n}:=\left( \pi _{1,n},...,\pi _{n,n}\right) $ and: $\pi _{m,n}=p\left( 
\frac{m}{n}\right) -p\left( \frac{m-1}{n}\right) $, $m=1,...,n.$ We note
that under our hypothesis, 
\begin{equation*}
\sum_{m=1}^{n}\pi _{m,n}=p\left( 1\right) -p\left( 0\right) =1.
\end{equation*}
Due to its asymmetry, the law of the biased $\mathbf{\nu }_{n}$ no longer is
exchangeable.

We now recall some well-known bias examples arising in population genetics. 
\newline

\emph{Example 3.1} (homographic model, selection). Assume 
\begin{equation}
p\left( x\right) =\left( 1+s\right) x/\left( 1+sx\right) ,  \label{se}
\end{equation}
where $s>-1$ is a selection parameter. This model arises when gene $A$
(respectively $a$), with frequency $x$ (respectively $1-x$), has fitness $%
1+s $ (respectively $1$). The case $s>0$ arises when gene of type $A$ is
selectively advantageous, whereas it is disadvantageous when $s\in \left(
-1,0\right) .$

\emph{Example 3.2}\textbf{\ }(selection with dominance). Assume 
\begin{equation}
p\left( x\right) =\frac{\left( 1+s\right) x^{2}+\left( 1+sh\right) x\left(
1-x\right) }{1+sx^{2}+2shx\left( 1-x\right) }.  \label{sd}
\end{equation}
In this model, genotype $AA$ (respectively $Aa$ and $aa$), with frequency $%
x^{2}$ (respectively $2x\left( 1-x\right) $ and $\left( 1-x\right) ^{2}$)
has fitness $1+s$ (respectively $1+sh$ and $1$). $h$ is a measure of the
degree of dominance of heterozygote $Aa$. We impose $s>-1$ and $sh>-1.$ Note
that the latter quantity can be put into the canonical form of deviation to
neutrality: 
\begin{equation*}
p\left( x\right) =x+sx\left( 1-x\right) \frac{h-x\left( 2h-1\right) }{%
1+sx^{2}+2shx\left( 1-x\right) }
\end{equation*}
where the ratio appearing in the right-hand-side is the ratio of the
difference of marginal fitnesses of $A$ and $a$ to their mean fitness. The
case $h=1/2$ corresponds to balancing selection with: $p\left( x\right) =x+%
\frac{s}{2}\frac{x\left( 1-x\right) }{1+sx}.$

\emph{Example 3.3}\textbf{\ }(quadratic model) With $c\in \left[ -1,1\right]
,$ a curvature parameter, one may choose: 
\begin{equation}
p\left( x\right) =x\left( 1+c-cx\right) .  \label{quad}
\end{equation}
If $c=1$, $p\left( x\right) =x\left( 2-x\right) =1-\left( 1-x\right) ^{2}$:
this bias appears in a discrete $2$-sex population model (M\"{o}hle, '$1994$%
, '$1998$)). We shall give below an interpretation of this quadratic model
when $c\in \left( 0,1\right] $ in terms of a joint one-way mutations and
neutrality effects model.

We can relax the assumption $p\left( 0\right) =0,$ $p\left( 1\right) =1$ by
assuming $0\leq p\left( 0\right) \leq $ $p\left( 1\right) \leq 1$, $p\left(
1\right) -p\left( 0\right) \in \left[ 0,1\right) .$

\emph{Example 3.4}\textbf{\ }(affine model) Take for example 
\begin{equation}
p\left( x\right) =\left( 1-\mu _{2}\right) x+\mu _{1}\left( 1-x\right) ,
\label{mu}
\end{equation}
where $\left( \mu _{1},\mu _{2}\right) $ are mutation probabilities,
satisfying $\mu _{1}\leq 1-\mu _{2}.$ It corresponds to the mutation scheme: 
$a\overset{}{\underset{\mu _{2}}{\overset{\mu _{1}}{\rightleftarrows }}}A$.
To avoid discussions of intermediate cases, we will assume that $p\left(
0\right) =\mu _{1}>0$ and $p\left( 1\right) <1$ ($\mu _{2}>0$)$.$ In this
case, the matrix $\Pi _{n}$ is irreducible and even primitive and all states
of this Markov chain are now recurrent. We have $\Bbb{P}\left( N_{r+1}>0\mid
N_{r}=0\right) =1-\left( 1-p\left( 0\right) \right) ^{n}>0$ and $\Bbb{P}%
\left( N_{r+1}<n\mid N_{r}=n\right) =1-p\left( 1\right) ^{n}>0$ and the
boundaries $\left\{ 0\right\} $ and $\left\{ n\right\} $ no longer are
strictly absorbing as there is a positive reflection probability inside the
domain $\left\{ 0,1,...,n\right\} .$\newline

For reasons to appear now, we shall be only interested in functions $q$ such
that $q\left( x\right) :=1-p\left( x\right) $ is a completely monotone
function (CM) on $\left( 0,1\right) $ that is, satisfying: 
\begin{equation*}
\left( -1\right) ^{l}q^{\left( l\right) }\left( x\right) \geq 0\text{, for
all }x\in \left( 0,1\right) ,
\end{equation*}
for all order-$l$ derivatives $q^{\left( l\right) }$ of $q$, $l\geq 0$. If $%
p\left( x\right) $ is such that $q$ is CM, we shall call it an admissible
bias mechanism.

\section{Non-neutral Wright-Fisher models and duality}

\emph{Preliminaries:} Let $\mathbf{v}_{n}:=\left( v\left( 0\right) ,v\left(
1\right) ,...,v\left( n\right) \right) $ be a $\left( n+1\right) -$vector of 
$\left[ 0,1\right] -$valued numbers. Define the backward difference operator 
$\nabla $ acting on $\mathbf{v}_{n}$ by: $\nabla v\left( m\right) =v\left(
m\right) -v\left( m-1\right) $, $m=1,...,n.$ We have $\nabla ^{0}v\left(
m\right) =v\left( m\right) ,$ $\nabla ^{2}v\left( m\right) =v\left( m\right)
-2v\left( m-1\right) +v\left( m-2\right) ,$ etc..., and, starting from the
endpoint $v\left( n\right) $ 
\begin{equation*}
\nabla ^{j}v\left( m\right) \mid _{m=n}=\sum_{l=0}^{j}\left( -1\right) ^{j-l}%
\binom{j}{l}v\left( n-l\right) \text{, }j=0,...,n.
\end{equation*}
Let $u$ be some continuous function: $\left[ 0,1\right] \rightarrow \left[
0,1\right] .$ Consider the $\left( n+1\right) -$vector $\mathbf{u}%
_{n}:=\left( u\left( \frac{0}{n}\right) ,...,u\left( \frac{m}{n}\right)
,...,u\left( \frac{n}{n}\right) \right) $, sampling $u$ at points $m/n.$ The
function $u$ is said to be $\nabla -$completely monotonic if $\left(
-1\right) ^{j}\nabla ^{j}u\left( \frac{m}{n}\right) \mid _{m=n}\geq 0$, for
all $j=0,...,n$ and all $n\geq 0.$ Let $\left( u^{1},u^{2}\right) $ be two
continuous functions on $\left[ 0,1\right] $. Let $u=u^{1}\cdot u^{2}$. With 
$\mathbf{u}_{n}$ the point-wise product of $\mathbf{u}_{n}^{1}$ and $\mathbf{%
u}_{n}^{2}$, assuming both functions $\left( u^{1},u^{2}\right) $ to be $%
\nabla -$completely monotonic, so will be $u,$ by the Leibniz rule. In
particular, if $u$ is $\nabla -$completely monotonic, so will be its
integral powers $u^{i},$ $i\in \Bbb{N}$. Our main result is:\newline

\textbf{Theorem:} \emph{Consider a non-neutral WF forward model }$\left(
N_{r};r\in \Bbb{N}\right) $\emph{\ on }$\left\{ 0,...,n\right\} $\emph{,
with continuous, non-decreasing bias }$p\left( x\right) ,$\emph{\
satisfying: } 
\begin{equation*}
0\leq p\left( 0\right) \leq p\left( 1\right) \leq 1,p\left( 1\right)
-p\left( 0\right) \in \left[ 0,1\right] .
\end{equation*}
\emph{This process has forward transition matrix: } 
\begin{equation*}
\Pi _{n}\left( k,k^{\prime }\right) =\Bbb{P}\left( \nu _{1,n}+...+\nu
_{k,n}=k^{\prime }\right) =\binom{n}{k^{\prime }}p\left( \frac{k}{n}\right)
^{k^{\prime }}\left( 1-p\left( \frac{k}{n}\right) \right) ^{n-k^{\prime }}.
\end{equation*}
\emph{There exists a Markov chain }$\left( A_{r};r\in \Bbb{N}\right) $\emph{%
\ on }$\left\{ 0,...,n\right\} $\emph{\ such that }$\left( N_{r},A_{r};r\in 
\Bbb{N}\right) $\emph{\ are dual with respect to }$\Phi _{n}^{2}\left(
m,k\right) =\binom{n-m}{k}/\binom{n}{k}$\emph{\ if and only if: }$%
x\rightarrow q\left( x\right) =1-p\left( x\right) $\emph{\ is completely
monotone on }$\left( 0,1\right) $\emph{. In this case, the transition
probability matrix of }$\left( A_{r};r\in \Bbb{N}\right) $\emph{\ is: } 
\begin{equation*}
P_{n}\left( i,j\right) =\binom{n}{j}\sum_{l=0}^{j}\left( -1\right) ^{j-l}%
\binom{j}{l}q\left( 1-\frac{l}{n}\right) ^{i}\geq 0.
\end{equation*}
$P_{n}$\emph{\ is a stochastic matrix if and only if }$p\left( 0\right) =0;$%
\emph{\ else, if }$p\left( 0\right) >0$\emph{, the matrix }$P_{n}$\emph{\ is
sub-stochastic.}\newline

\textbf{Proof:} Developing $\left[ \Phi _{n}^{2}\right] ^{-1}\Pi _{n}\Phi
_{n}^{2}=P_{n}^{\prime },$ we easily obtain: 
\begin{eqnarray*}
P_{n}^{^{\prime }}\left( j,i\right) &=&P_{n}\left( i,j\right) =\binom{n}{j}%
\sum_{l=0}^{j}\left( -1\right) ^{j-l}\binom{j}{l}\left[ 1-p\left( \frac{n-l}{%
n}\right) \right] ^{i} \\
&=&\binom{n}{j}\left( -1\right) ^{j}\nabla ^{j}\left( q\left( \frac{m}{n}%
\right) ^{i}\right) \mid _{m=n}
\end{eqnarray*}
This entry is non-negative if and only if $\left( -1\right) ^{j}\nabla
^{j}\left( q\left( \frac{m}{n}\right) ^{i}\right) \mid _{m=n}\geq 0$, for
all $i,j=0,...,n.$ But, due to the above argument on $\nabla -$complete
monotonicity of integral powers, this will be the case if and only if $%
\left( -1\right) ^{j}\nabla ^{j}\left( q\left( \frac{m}{n}\right) \right)
\mid _{m=n}\geq 0$, for all $j=0,...,n$. As this must be true for arbitrary
value of population size $n$, function $q$ has to be $\nabla -$completely
monotonic. Adapting now the arguments of Theorem $2$ developed in Feller '$%
1971$, page $223$, for absolutely monotone functions on $\left( 0,1\right) $%
, this will be the case if and only if $x\rightarrow q\left( x\right)
:=1-p\left( x\right) $ is a completely monotone function on $\left(
0,1\right) $ in the sense that: 
\begin{equation*}
\left( -1\right) ^{l}q^{\left( l\right) }\left( x\right) \geq 0\text{, for
all }x\in \left( 0,1\right) \text{, }l\in \Bbb{N}.
\end{equation*}

Next, since $\left( I-\nabla \right) u\left( m\right) =u\left( m-1\right) $
is a simple back-shift, 
\begin{equation*}
\sum_{j=0}^{n}P_{n}^{^{\prime }}\left( j,i\right) =\sum_{j=0}^{n}P_{n}\left(
i,j\right) =\left( I-\nabla \right) ^{n}\left( q\left( \frac{m}{n}\right)
^{i}\right) \mid _{m=n}=q\left( 0\right) ^{i}
\end{equation*}
and, if $q$ is CM, $P_{n}$ is a stochastic matrix if and only if $q\left(
0\right) =1;$ else, if $q\left( 0\right) <1$, the matrix $P_{n}$ is
sub-stochastic.

We note that the first column of the matrix $P_{n}$ is $P_{n}\left(
i,0\right) =q\left( 1\right) ^{i}$ whereas its first line is: $P_{n}\left(
0,j\right) =\delta _{0,j}$, expressing, as required, that the state $0$ of $%
\left( A_{r};r\in \Bbb{N}\right) $ is absorbing. $\square $

\section{Examples}

We show here that most of the simplest non-neutrality mechanisms used in the
literature fall within the class which we would like to draw the attention
on, or are amenable to it via some `reciprocal transformation' which we
define. Elementary algebraic manipulations on CM functions allows to exhibit
a vast class of unsuspected mechanisms. Note that in some cases, their
biological relevance remains to be elucidated. The results presented in this
Section seem to be new. They serve as an illustration of our theorem.

\subsection{Elementary mechanisms}

Assume first $p\left( x\right) =x$ corresponding to the simple neutral case.
Then $q\left( x\right) =1-x$ is completely monotone on $\left( 0,1\right) $.
With $S_{i,j}$ the second kind Stirling numbers, we get a lower left
triangular stochastic transition matrix 
\begin{eqnarray*}
P_{n}\left( i,j\right) &=&\binom{n}{j}\sum_{l=0}^{j}\left( -1\right) ^{j-l}%
\binom{j}{l}\left( \frac{l}{n}\right) ^{i}=\left( n\right) _{j}\cdot
n^{-i}\cdot S_{i,j}\text{, }j\leq i \\
P_{n}\left( i,j\right) &=&0\text{, else.}
\end{eqnarray*}
The diagonal terms (eigenvalues) are all distinct with $P_{n}\left(
i,i\right) =\left( n\right) _{i}\cdot n^{-i}.$ The matrix $P_{n}$ is
stochastic. Due to triangularity, ancestral process is a pure death Markov
process which may be viewed as a discrete coalescence tree. \newline

From example 3.4 (mutation). Assume (\ref{mu}) holds: $p\left( x\right)
=\left( 1-\mu _{2}\right) x+\mu _{1}\left( 1-x\right) $ where $\left( \mu
_{1},\mu _{2}\right) $ are mutation probabilities. Then, with $\kappa
:=1-\left( \mu _{1}+\mu _{2}\right) $, $q\left( x\right) =1-\mu _{1}-$ $%
\kappa x$ is completely monotone on $\left( 0,1\right) $ if and only $\mu
_{1}\leq 1-\mu _{2}$ ($\kappa \geq 0$)$.$ In this case, $P_{n}$ is again
lower left triangular (a pure death process). We have 
\begin{equation}
P_{n}\left( i,j\right) =\binom{n}{j}\sum_{l=0}^{j}\left( -1\right) ^{j-l}%
\binom{j}{l}\left( \mu _{2}+\kappa \frac{l}{n}\right) ^{i}=:\left( n\right)
_{j}\cdot n^{-i}\cdot S_{i,j}^{\mu _{2}}\left( \kappa /n\right) \text{, }%
j\leq i  \label{mut}
\end{equation}
\begin{equation*}
P_{n}\left( i,j\right) =0,\text{ else,}
\end{equation*}
in terms of generalized Stirling numbers $S_{i,j}^{\mu _{2}}\left( \kappa
/n\right) $. We have $P_{n}\left( i,i\right) =\left( n\right) _{i}\left( 
\frac{\kappa }{n}\right) ^{i}$ and the spectrum of $P_{n}$ is real. When $%
\mu _{1}>0$, this matrix is sub-stochastic with $\sum_{j=0}^{n}P_{n}\left(
i,j\right) =\left( 1-\mu _{1}\right) ^{i}$.

A particular case deals with one-way mutations ($\mu _{1}+\mu _{2}>0$, $\mu
_{1}\cdot \mu _{2}=0$):

If $\mu _{2}=0,$ $P_{n}\left( i,j\right) =\left( 1-\mu _{1}\right) ^{i}\cdot
\left( n\right) _{j}\cdot n^{-i}\cdot S_{i,j}$, $j\leq i$, $=0$, else.
Further, $\sum_{j=0}^{n}P_{n}\left( i,j\right) =\left( 1-\mu _{1}\right)
^{i}<1.$

If $\mu _{1}=0,$ $P_{n}\left( i,j\right) =\left( n\right) _{j}\cdot
n^{-i}\cdot S_{i,j}^{\mu _{2}}\left( \left( 1-\mu _{2}\right) /n\right) $, $%
j\leq i$, $=0$, else. The corresponding matrix $P_{n}$ is stochastic.\newline

From example 3.3 (quadratic). Assume $p\left( x\right) =x\left(
1+c-cx\right) $, as in (\ref{quad}). Then $q\left( x\right) =\left(
1-x\right) \left( 1-cx\right) $ which is completely monotone if and only if $%
c\in \left[ 0,1\right] $. The case $c=0$ is the neutral case, whereas $c=1$
appears in a $2$-sex model of M\"{o}hle. In this quadratic case, since $%
\nabla ^{j}\left( q\left( \frac{m}{n}\right) ^{i}\right) =0$ if $j>2i,$ then 
$P_{n}\left( i,j\right) =0$ if $j>2i$ and so $P_{n}$ is a Hessenberg-like
matrix$.$ Note that $\sum_{j=0}^{n}P_{n}\left( i,j\right) =q\left( 0\right)
^{i}=1$.\newline

From the selection example 3.1, when (\ref{se}) holds 
\begin{equation*}
p\left( x\right) =\left( 1+s\right) x/\left( 1+sx\right) ,
\end{equation*}
$q\left( x\right) =1-p\left( x\right) =\left( 1-x\right) /\left( 1+sx\right) 
$ is CM whenever selection parameter $s>0$. The induced matrix $P_{n}$ is
stochastic. It is no longer lower left triangular so that the ancestral no
longer is a pure death process, rather a birth and death process. The
induced coalescence pattern no longer is a discrete tree, but rather a graph
(a discrete version of the ancestral selection graph of Neuhauser-Krone '$%
1997$).\newline

From example 3.2 (selection with dominance). The corresponding mechanism (%
\ref{sd}) with parameters ($s,h$) satisfying $s>-1$ and $sh>-1$ is CM if and
only if $s>0$ and $h\in \left( 0,1/2\right) $. The case $h\in \left(
0,1\right) $ corresponds to directional selection where genotype $AA$ has
highest fitness compared to $aa$'s and the heterozygote class $Aa$ has
intermediate fitness compared to both homozygote classes. In this situation,
marginal fitness of $A$ exceeds the one of $a$ and selective sweep is
expected. When $h\in \left( 0,1/2\right) $, allele $A$ is dominant to $a$,
whereas when $h\in \left( 1/2,1\right) $, allele $A$ is recessive to $a$ (a
stabilizing effect slowing down the sweep). Critical value $h=1/2$ is a case
of pure genic balancing selection.\newline

\emph{Example 5.1.1} Consider the mechanism $p\left( x\right) =x^{\gamma }$
for some $\gamma >0$. The function $q\left( x\right) =1-p\left( x\right) $
is CM if and only if $\gamma \in \left( 0,1\right) .$ Although this model
seems quite appealing, we could find no reference to it in the specialized
mathematical genetics' literature.

\emph{Example 5.1.2}\textbf{\ }(Reciprocal mechanism) If $p\left( x\right) $
is not admissible in that $q$ is not CM, it can be that $\overline{p}\left(
x\right) :=1-p\left( 1-x\right) $ is itself admissible. As observed before,
if $N_{r}\left( m\right) $ has transition probabilities given by Bin$\left(
n,p\left( k/n\right) \right) ,$ $\overline{p}\left( x\right) $ arises in the
transition probabilities of $\overline{N}_{r}\left( m\right) :=n-N_{r}\left(
m\right) .$ Indeed, such transitions are Bin$\left( n,\overline{p}\left(
k/n\right) \right) $ distributed.

If $p\left( x\right) $ is the selection mechanism of example 3.1, (\ref{se}%
), with $s\in \left( -1,0\right) $ (not admissible), $\overline{p}\left(
x\right) =\left( 1+\overline{s}\right) x/\left( 1+\overline{s}x\right) $ is
itself an admissible selection mechanism because it has reciprocal selection
parameter $\overline{s}=-s/\left( 1+s\right) >0.$ If $p\left( x\right) $ is
the mechanism of example 3.2, namely (\ref{sd}), with parameters $\left(
s,h\right) $, then $\overline{p}\left( x\right) $ is itself a selection with
dominance mechanism with reciprocal parameters $\overline{s}=-s/\left(
1+s\right) $ and $\overline{h}=1-h.$ Assuming $\left( s<0,h\in \left(
1/2,1\right) \right) $, $p\left( x\right) $ is not admissible whereas $%
\overline{p}\left( x\right) $ is because $\overline{s}>0$ and $\overline{h}%
\in \left( 0,1/2\right) $. Similarly, when $\gamma \in \left( 0,1\right) $,
the mechanism $p\left( x\right) =1-\left( 1-x\right) ^{\gamma }$ is not
admissible but, from example 5.1.1, $\overline{p}\left( x\right) :=1-p\left(
1-x\right) =x^{\gamma }$ is.

\subsection{Bias mechanisms with mutational effects}

Let $p_{M}\left( x\right) =\left( 1-\mu _{2}\right) x+\mu _{1}\left(
1-x\right) $ be the mutational bias mechanism (with $\kappa =1-\left( \mu
_{1}+\mu _{2}\right) \geq 0$). Let $p\left( x\right) $ be a bias mechanism
such that $q\left( x\right) $ is CM with $p\left( 1\right) -p\left( 0\right)
=1$. Then 
\begin{equation*}
\widetilde{p}_{M}\left( x\right) =p_{M}\left( p\left( x\right) \right)
\end{equation*}
is such that $\widetilde{q}_{M}\left( x\right) :=1-\widetilde{p}_{M}\left(
x\right) $ is CM. It is therefore admissible and adds mutational effects to
the primary mechanism $p\left( x\right) .$ For example, 
\begin{equation*}
\widetilde{p}_{M}\left( x\right) =\frac{\mu _{1}+x\left( \left( 1+s\right)
\left( 1-\mu _{2}\right) -\mu _{1}\right) }{1+sx}
\end{equation*}
is a mechanism of selection combined with mutational effects. We have $%
\widetilde{p}_{M}\left( 0\right) =\mu _{1}$, $\widetilde{p}_{M}\left(
1\right) =1-\mu _{2}.$ The mechanisms $\widetilde{p}_{M}\left( x\right) $
obtained in this way all share the specificity: $\widetilde{p}_{M}\left(
1\right) -$ $\widetilde{p}_{M}\left( 0\right) =:\kappa <1.$

Note that, except for the mutational affine mechanism, it is not true in
general that whenever $p^{1}\left( x\right) $ and $p^{2}\left( x\right) $
are two admissible bias mechanisms, then $p^{1}\left( p^{2}\left( x\right)
\right) $ is admissible.

\subsection{Joint bias effects and Compound bias}

Let $p^{1}\left( x\right) $ and $p^{2}\left( x\right) $ be two admissible
bias in that $q^{1}\left( x\right) :=1-p^{1}\left( x\right) $ and $%
q^{2}\left( x\right) :=1-p^{2}\left( x\right) $ are both completely
monotone. Then 
\begin{equation*}
q\left( x\right) =q^{1}\left( x\right) q^{2}\left( x\right) \text{ is CM.}
\end{equation*}
Thus, with $x_{1}\circ x_{2}:=x_{1}+x_{2}-x_{1}x_{2},$ the probabilistic
product in $\left[ 0,1\right] $%
\begin{equation*}
\left( p^{1}\left( x\right) ,p^{2}\left( x\right) \right) \rightarrow
p\left( x\right) =p^{1}\left( x\right) \circ p^{2}\left( x\right) .
\end{equation*}
Whenever a WF model is considered with bias $p\left( x\right) =p^{1}\left(
x\right) \circ p^{2}\left( x\right) $ obtained from two distinct bias $%
p^{1}\left( x\right) $ and $p^{2}\left( x\right) $, we call it a WF model
with \emph{joint bias effect}.\newline

\emph{Example 5.3.1} (Joint selection and mutational effects). Let $%
p^{1}\left( x\right) =p_{M}\left( x\right) $ and $p^{2}\left( x\right)
=\left( 1+s\right) x/\left( 1+sx\right) $. We get 
\begin{equation*}
q\left( x\right) =\frac{\left( 1-\mu _{1}-\kappa x\right) \left( 1-x\right) 
}{1+sx}\text{ and }p\left( x\right) =\frac{\mu _{1}+x\left( s+1-\mu
_{1}+\kappa \right) -\kappa x^{2}}{1+sx},
\end{equation*}
with $p\left( 0\right) =\mu _{1}$, $p\left( 1\right) =1.$ This mechanism
differs from the traditional mechanism of selection combined with mutational
effects.

\emph{Example 5.3.2} (Joint mutation and neutral effects). Let $p^{1}\left(
x\right) =\left( 1-\mu _{2}\right) x+\mu _{1}\left( 1-x\right) $ and $%
p^{2}\left( x\right) =x$. We get 
\begin{equation*}
q\left( x\right) =\left( 1-x\right) \left( 1-\mu _{1}-\kappa x\right) \text{
and }p\left( x\right) =\mu _{1}+x\left( 1-\mu _{1}+\kappa \left( 1-x\right)
\right) ,
\end{equation*}
with $p\left( 0\right) =\mu _{1}$, $p\left( 1\right) =1.$ When $\mu _{1}=0$
(one-way mutations), we recover the quadratic mechanism with curvature
parameter $c=1-\mu _{2}.$ This finding justifies some interest into the
quadratic mechanisms with $c\neq 1$.

With $j=1,2$, the reproduction law of each elementary effect is $\mathbf{\nu 
}_{n}^{j}\overset{d}{\sim }$ Multin$\left( n;\mathbf{\pi }_{n}^{j}\right) ,$
where $\pi _{m,n}^{j}=p^{j}\left( \frac{m}{n}\right) -p^{j}\left( \frac{m-1}{%
n}\right) $, $m=1,...,n.$ Then, $\mathbf{\nu }_{n}\overset{d}{\sim }$ Multin$%
\left( n;\mathbf{\pi }_{n}\right) ,$ $\pi _{m,n}=p\left( \frac{m}{n}\right)
-p\left( \frac{m-1}{n}\right) $, $m=1,...,n$, where $\mathbf{\pi }_{n}:=%
\mathbf{\pi }_{n}^{1}\odot \mathbf{\pi }_{n}^{2}$ is easily obtained
component-wise by: 
\begin{equation*}
\pi _{m,n}=\pi _{m,n}^{1}\sum_{l=1}^{m}\pi _{l,n}^{2}+\pi
_{m,n}^{2}\sum_{l=1}^{m}\pi _{l,n}^{1},m=1,...,n.
\end{equation*}
We let: $\mathbf{\nu }_{n}:=\mathbf{\nu }_{n}^{1}\odot \mathbf{\nu }_{n}^{2}%
\overset{d}{\sim }$ Multin$\left( n;\mathbf{\pi }_{n}^{1}\odot \mathbf{\pi }%
_{n}^{2}\right) $. It is the reproduction law of a WF model obtained jointly
from the two bias $p^{1}\left( x\right) $ and $p^{2}\left( x\right) $.%
\newline

Let $\phi \left( x\right) :$ $\left( 0,1\right) \rightarrow \left(
0,1\right) $ be an absolutely monotone function satisfying: $\phi ^{\left(
l\right) }\left( x\right) \geq 0$ for all $l-$th derivatives $\phi ^{\left(
l\right) }$ of $\phi $, all $x\in \left( 0,1\right) $. Such functions are
well-known to be probability generating functions (pgfs) of $\Bbb{N}-$valued
random variables, say $N$, that is to say: $\phi \left( x\right) =\Bbb{E}%
\left[ x^{N}\right] .$ Clearly, if $q$ is CM on $\left( 0,1\right) $, then
so is: $q_{\phi }\left( x\right) :=\phi \left( q\left( x\right) \right) .$
Thus $p_{\phi }\left( x\right) :=1-\phi \left( 1-p\left( x\right) \right) $
is an admissible bias mechanism in that $q_{\phi }\left( x\right)
:=1-p_{\phi }\left( x\right) $ is CM. We call it a \emph{compound} bias.%
\newline

\emph{Example 5.3.3} The general mechanism with mutational effects is in
this class. Indeed, 
\begin{equation*}
\widetilde{q}_{M}\left( x\right) =1-\widetilde{p}_{M}\left( x\right)
=1-p_{M}\left( p\left( x\right) \right) =1-p_{M}\left( 1-q\left( x\right)
\right)
\end{equation*}
and so $\phi \left( x\right) =1-p_{M}\left( 1-x\right) =1-\left( 1-\mu
_{2}\right) \left( 1-x\right) -\mu _{1}x=\mu _{2}+\kappa x$ which is
absolutely monotone as soon as $\kappa =1-\left( \mu _{1}+\mu _{2}\right)
\geq 0.$

\emph{Example 5.3.4} With $\theta >0$, taking $\phi \left( x\right)
=e^{-\theta \left( 1-x\right) }$ or $\left( e^{\theta x}-1\right) /\left(
e^{\theta }-1\right) $, the pgf of a Poisson (or shifted-Poisson) random
variable, $p_{\phi }\left( x\right) =1-\phi \left( 1-p\left( x\right)
\right) $ is admissible if $p\left( x\right) $ is. Note that if $q$ is of
the form $q_{\phi }$ where $\phi $ is the pgf of a Poisson random variable,
then $q_{\phi }\left( x\right) ^{\alpha }$ is admissible for all $\alpha >0$%
, a property reminiscent of infinite divisibility for pgfs. Taking $\phi
\left( x\right) =\left( 1-\pi \right) /\left( 1-\pi x\right) $ or $x\left(
1-\pi \right) /\left( 1-\pi x\right) $, $\pi \in \left( 0,1\right) $, the
pgf of a geometric (or shifted-geometric) random variable, $p_{\phi }\left(
x\right) =sp\left( x\right) /\left( 1+sp\left( x\right) \right) $ or $\left(
s+1\right) p\left( x\right) /\left( 1+sp\left( x\right) \right) $ is
admissible if $p\left( x\right) $ is (with $s=\pi /\left( 1-\pi \right) >0$%
). In the external latter mechanism, one recognizes the one in (\ref{se})
occurring in the model with selection of example 3.1.

\emph{Example 5.3.5} Let $p\left( x\right) =x^{\gamma }$ with $\gamma \in
\left( 0,1\right) $ as in example 5.1.1$.$ Then $p_{\phi }\left( x\right)
=1-q_{\phi }\left( x\right) $ where $q_{\phi }\left( x\right) =e^{-\theta
\left( 1-q\left( x\right) \right) }=e^{-\theta x^{\gamma }},$ $\theta >0$,
is admissible. Note that $p_{\phi }\left( x\right) \underset{x\downarrow 0}{%
\sim }\theta x^{\gamma }$. The reciprocal function $\overline{p}_{\phi
}\left( x\right) =q_{\phi }\left( 1-x\right) =e^{-\theta \left( 1-x\right)
^{\gamma }}$ also interprets as an absolutely monotone discrete-stable pgf
(see Steutel, van Harn, '$1979$). It is not admissible.

Proceeding in this way, one can produce a wealth of admissible bias
probabilities $p_{\phi }$, the signification of which in Population Genetics
remaining though to be pinpointed, in each specific case study.

\section{Limit laws}

Consider a WF model $\left( N_{r};r\in \Bbb{N}\right) $ on $\left\{
0,...,n\right\} $ with forward transition matrix: 
\begin{equation*}
\Pi _{n}\left( k,k^{\prime }\right) =\binom{n}{k^{\prime }}p\left( \frac{k}{n%
}\right) ^{k^{\prime }}\left( 1-p\left( \frac{k}{n}\right) \right)
^{n-k^{\prime }},
\end{equation*}
with admissible bias $p\left( x\right) .$ Define $\left( A_{r};r\in \Bbb{N}%
\right) $ as the dual Markov chain on $\left\{ 0,...,n\right\} $ with
transition probability: 
\begin{equation*}
P_{n}\left( i,j\right) =\binom{n}{j}\sum_{l=0}^{j}\left( -1\right) ^{j-l}%
\binom{j}{l}q\left( 1-\frac{l}{n}\right) ^{i}.
\end{equation*}
Then, $\left( N_{r},A_{r};r\in \Bbb{N}\right) $ are dual with respect to $%
\Phi _{n}\left( m,k\right) :=\Phi _{n}^{2}\left( m,k\right) =\binom{n-m}{k}/%
\binom{n}{k}$, to wit: 
\begin{equation*}
\text{ }\Bbb{E}_{m}\left[ \binom{n-N_{r}}{k}/\binom{n}{k}\right] =\Bbb{E}%
_{k}\left[ \binom{n-m}{A_{r}}/\binom{n}{A_{r}}\right] =\Bbb{E}_{k}\left[ 
\binom{n-A_{r}}{m}/\binom{n}{m}\right] .
\end{equation*}
We shall distinguish two cases.\newline

\textbf{Case 1:} Assume first that 
\begin{equation*}
N_{r}\overset{d}{\rightarrow }N_{\infty }\text{ as }r\uparrow \infty ,\text{
independently of }N_{0}=m\geq 1.
\end{equation*}
Let $\pi _{\infty }\left( i\right) =\Bbb{P}\left( N_{\infty }=i\right) $ and 
$\mathbf{\pi }_{\infty }=\left( \pi _{\infty }\left( 0\right) ,...,\pi
_{\infty }\left( n\right) \right) ^{\prime }.$ The line vector $\mathbf{\pi }%
_{\infty }^{\prime }$ is the left eigenvector of $\Pi _{n}$ associated to
the eigenvalue $1:$ $\mathbf{\pi }_{\infty }^{\prime }=\mathbf{\pi }_{\infty
}^{\prime }\Pi _{n}.$ It is the (unique) invariant probability measure
(stationary distribution) of $\left( N_{r};r\in \Bbb{N}\right) .$

If this stationary distribution exists, then, using duality formula,
necessarily, $A_{r}\rightarrow 0$ as $r\uparrow \infty $ with probability $%
\Bbb{P}_{k}\left( A_{\infty }=0\right) =:\rho _{\infty }\left( k\right) <1.$
The numbers $\rho _{\infty }\left( k\right) $ are the extinction
probabilities of the dual process started at $k.$ As is well-known, $\mathbf{%
\rho }_{\infty }=\left( \rho _{\infty }\left( 0\right) ,...,\rho _{\infty
}\left( n\right) \right) ^{\prime }$ is the unique solution to $\left(
I-P_{n}\right) \mathbf{\rho }_{\infty }=0$ with $\rho _{\infty }\left(
0\right) =1.$ \newline

\textbf{Remark:} Typical situations where $\left( N_{r};r\in \Bbb{N}\right) $
has an invariant measure is when mutational effects are present, and more
generally when the bias mechanism satisfies $p\left( 0\right) >0$ and $%
p\left( 1\right) <1$. In this situation, the forward stochastic transition
matrix $\Pi _{n}$ has an algebraically simple dominant eigenvalue $1$. By
Perron-Frobenius theorem: 
\begin{equation*}
\lim_{r\uparrow \infty }\Pi _{n}^{r}=\mathbf{1\pi }_{\infty }^{\prime },
\end{equation*}
where $\mathbf{1}^{\prime }=\left( 1,...,1\right) .$ The invariant
probability measure can be approximated by subsequent iterates of $\Pi _{n}$%
, the convergence being exponentially fast, with rate governed by the second
largest eigenvalue. Of course, detailed balance (stating that $\pi _{k}\Pi
_{n}\left( k,k^{\prime }\right) =\pi _{k^{\prime }}\Pi _{n}\left( k^{\prime
},k\right) $) does not hold here and the forward chain in equilibrium is not
time-reversible.

In these recurrent cases, the dual ancestral process $A_{r}$ started at $k$
gets extinct with probability $\rho _{\infty }\left( k\right) $. The numbers 
$1-\rho _{\infty }\left( k\right) $ are the probabilities that it gets
killed before getting extinct; in other words, $1-\rho _{\infty }\left(
k\right) $ are the probabilities that $A_{r}$ first hits an extra coffin
state, say $\left\{ \partial \right\} ,$ before hitting $\left\{ 0\right\} $%
. $\square $\newline

In terms of moments, by the duality formula, we conclude that: 
\begin{equation*}
\binom{n}{k}^{-1}\Bbb{E}\left[ \binom{n-N_{\infty }}{k}\right] =\rho
_{\infty }\left( k\right) =\Bbb{P}_{k}\left( A_{\infty }=0\right) ,
\end{equation*}
relating $k-$factorial moments of $n-N_{\infty }$ to the extinction
probabilities of $A_{r}$ given $A_{0}=k.$ We also have 
\begin{equation*}
\sum_{k=0}^{n}v^{k}\Bbb{E}\left[ \binom{n-N_{\infty }}{k}\right] =\Bbb{E}%
\left[ \left( 1+v\right) ^{n-N_{\infty }}\right] =\sum_{k=0}^{n}\binom{n}{k}%
\rho _{\infty }\left( k\right) v^{k}
\end{equation*}
and so, the probability generating function of $N_{\infty }$ can be
expressed as ($u\in \left[ 0,1\right] $): 
\begin{equation*}
\Bbb{E}\left[ u^{N_{\infty }}\right] =\sum_{k=0}^{n}\binom{n}{k}\rho
_{\infty }\left( k\right) u^{n-k}\left( 1-u\right) ^{k},
\end{equation*}
in terms of the Bernstein-B\'{e}zier polynomial of $\left( \rho _{\infty
}\left( n-k\right) ;\text{ }k=0,...,n\right) .$\newline

Let $\mathbf{\rho }_{\infty }=\left( \rho _{\infty }\left( 0\right)
,...,\rho _{\infty }\left( n\right) \right) ^{\prime }.$ The vector $\mathbf{%
\rho }_{\infty }$ is the right eigenvector of $P_{n}$ associated to the
eigenvalue $1:$ $\mathbf{\rho }_{\infty }=P_{n}\mathbf{\rho }_{\infty }.$ In
this case, the matrix $P_{n}$ is sub-stochastic and the extinction
probability of $\left( A_{r};r\in \Bbb{N}\right) $ given $A_{0}=k$ is less
than one. Thanks to duality, we have: 
\begin{equation*}
\Pi _{n}\Phi _{n}=\Phi _{n}P_{n}^{\prime }.
\end{equation*}
where the matrix $\Phi _{n}$ is symmetric whereas the matrix $\Phi _{n}^{-1}$
is symmetric right-lower triangular, with: 
\begin{equation*}
\Phi _{n}\left( m,k\right) =\binom{n-m}{k}/\binom{n}{k}=\binom{n-k}{m}/%
\binom{n}{m}
\end{equation*}
\begin{equation*}
\text{ }\Phi _{n}^{-1}\left( i,j\right) =\left( -1\right) ^{i+j-n}\binom{i}{%
n-j}\binom{n}{i}=\left( -1\right) ^{i+j-n}\binom{j}{n-i}\binom{n}{j}.
\end{equation*}
Thus, 
\begin{equation*}
\mathbf{\pi }_{\infty }^{\prime }\Pi _{n}\Phi _{n}=\mathbf{\pi }_{\infty
}^{\prime }\Phi _{n}=\mathbf{\pi }_{\infty }^{\prime }\Phi _{n}P_{n}^{\prime
},
\end{equation*}
showing that $\mathbf{\rho }_{\infty }$ and $\mathbf{\pi }_{\infty }$ are
related through: 
\begin{equation*}
\mathbf{\rho }_{\infty }=\Phi _{n}\mathbf{\pi }_{\infty }\text{ or }\mathbf{%
\pi }_{\infty }=\Phi _{n}^{-1}\mathbf{\rho }_{\infty }.
\end{equation*}
\emph{The knowledge of the invariant measure }$\mathbf{\pi }_{\infty }$\emph{%
\ of the forward process allows one to compute the extinction probabilities }%
$\mathbf{\rho }_{\infty }$\emph{\ of the dual backward ancestral process
(and conversely).}\newline

\emph{Example:} Consider the discrete WF model with mutations of example
3.4. In this case, $N_{r}\overset{d}{\rightarrow }N_{\infty }$ as $r\uparrow
\infty ,$ regardless of $N_{0}=m$ and $\left( N_{r};r\in \Bbb{N}\right) $
has an invariant measure which is difficult to compute. Looking at the
backward process, the matrix $P_{n}$ is sub-stochastic (if $\mu _{1}>0$) and
lower-left triangular. Due to triangularity, the right eigenvector $\mathbf{%
\rho }_{\infty }$ of $P_{n}$ can easily be computed explicitly in terms of $%
\left( P_{n}\left( i,j\right) \text{; }j\leq i\right) ,i=0,...,n.$ We
therefore get the following alternating expression for the invariant
measure: 
\begin{equation*}
\pi _{\infty }\left( i\right) =\binom{n}{i}\sum_{j=0}^{i}\left( -1\right)
^{i-j}\binom{i}{j}\rho _{\infty }\left( n-j\right) .
\end{equation*}
Concerning moments, for instance, we have $\rho _{\infty }\left( 1\right)
=\mu _{2}/\left( \mu _{1}+\mu _{2}\right) $ so that $\Bbb{E}\left[ N_{\infty
}\right] =n\mu _{1}/\left( \mu _{1}+\mu _{2}\right) $; from (\ref{mut}) we
also have: 
\begin{equation*}
\rho _{\infty }\left( 2\right) =\frac{\mu _{2}\left[ n\mu _{2}\left(
1+\kappa \right) +\kappa ^{2}\right] }{\left[ \left( 1-\kappa \right) \left(
n-\left( n-1\right) \kappa ^{2}\right) \right] }=\frac{1}{n\left( n-1\right) 
}\left( n\left( n-1\right) -n\frac{\left( 2n-1\right) \mu _{1}}{\mu _{1}+\mu
_{2}}+\Bbb{E}\left[ N_{\infty }^{2}\right] \right)
\end{equation*}
allowing to compute $\Bbb{E}\left[ N_{\infty }^{2}\right] $ and then the
variance of $N_{\infty }.$ We get: 
\begin{equation*}
\sigma ^{2}\left( N_{\infty }\right) =\frac{n^{2}\mu _{1}\mu _{2}}{\left(
\mu _{1}+\mu _{2}\right) ^{2}\left( 2n\left( \mu _{1}+\mu _{2}\right)
+1\right) }+o\left( n\right) \underset{n\uparrow \infty }{\sim }\frac{\mu
_{1}\mu _{2}}{2\left( \mu _{1}+\mu _{2}\right) ^{3}}n,
\end{equation*}
suggesting (when $\mu _{1}\mu _{2}>0$) a Central Limit Theorem for $%
N_{\infty }$ as $n$ grows large: 
\begin{equation*}
\frac{1}{\sqrt{n}}\left( N_{\infty }-n\frac{\mu _{1}}{\mu _{1}+\mu _{2}}%
\right) \overset{d}{\underset{n\uparrow \infty }{\rightarrow }}\mathcal{N}%
\left( 0,\frac{\mu _{1}\mu _{2}}{2\left( \mu _{1}+\mu _{2}\right) ^{3}}%
\right) .\text{ }
\end{equation*}
\newline

\textbf{Case 2.} Conversely, assume now that given $N_{0}=m$%
\begin{equation*}
N_{r}\overset{d}{\rightarrow }0\text{ as }r\uparrow \infty ,\text{ with
probability }\Bbb{P}_{m}\left( N_{\infty }=0\right) =:\rho _{\infty }\left(
m\right) ,
\end{equation*}
so that boundaries $\left\{ 0,n\right\} $ are absorbing. Then, the ancestral
process $\left( A_{r};r\in \Bbb{N}\right) $ possesses an invariant
distribution, in that: 
\begin{equation*}
A_{r}\overset{d}{\rightarrow }A_{\infty }\text{ as }r\uparrow \infty ,\text{
independently of }A_{0}=k\in \left[ n\right] .
\end{equation*}
In terms of moments, the duality formula means that: 
\begin{equation*}
\binom{n}{m}^{-1}\Bbb{E}\left[ \binom{n-A_{\infty }}{m}\right] =\rho
_{\infty }\left( m\right) =\Bbb{P}_{m}\left( N_{\infty }=0\right) ,
\end{equation*}
relating $m-$factorial moments of $n-A_{\infty }$ to the extinction
probabilities of $N_{r}$ given $N_{0}=m.$ Stated differently, the
probability generating function of $A_{\infty }$ is ($u\in \left[ 0,1\right] 
$): 
\begin{equation*}
\Bbb{E}\left[ u^{A_{\infty }}\right] =\sum_{m=0}^{n}\binom{n}{m}\rho
_{\infty }\left( m\right) u^{n-m}\left( 1-u\right) ^{m}.
\end{equation*}
\newline

Let $\pi _{\infty }\left( i\right) =\Bbb{P}\left( A_{\infty }=i\right) $,
with $\mathbf{\pi }_{\infty }^{\prime }=\mathbf{\pi }_{\infty }^{\prime
}P_{n}.$ Then, using duality$,$ $\mathbf{\rho }_{\infty }$ is the right
eigenvector of $\Pi _{n}$ associated to the eigenvalue $1:$ $\mathbf{\rho }%
_{\infty }=\Pi _{n}\mathbf{\rho }_{\infty }.$ Thus, $\mathbf{\rho }_{\infty
} $ and $\mathbf{\pi }_{\infty }$ are related through: 
\begin{equation*}
\mathbf{\rho }_{\infty }=\Phi _{n}\mathbf{\pi }_{\infty }\text{ or }\mathbf{%
\pi }_{\infty }=\Phi _{n}^{-1}\mathbf{\rho }_{\infty }.
\end{equation*}

\emph{The knowledge of the extinction probabilities }$\mathbf{\rho }_{\infty
}$\emph{\ of the forward process allows one to compute the invariant measure 
}$\mathbf{\pi }_{\infty }$ \emph{of the dual backward ancestral process (and
conversely).}\newline

\emph{Examples}: Typical situations where boundaries $\left\{ 0,n\right\} $
are absorbing to $\left( N_{r};r\in \Bbb{N}\right) $ occur when $p\left(
0\right) =0$ and $p\left( 1\right) =1$. The simplest case is the neutral
case, but the non-neutral selection and selection with dominance mechanisms
or the quadratic mechanism (examples 3.1, 3.2 and 3.3) are also in this
class. For instance:

$\left( i\right) $ In the neutral case, $\rho _{\infty }\left( m\right)
=1-m/n.$ Thus, $\pi _{\infty }\left( i\right) =\binom{n}{i}%
\sum_{j=0}^{i}\left( -1\right) ^{i-j}\binom{i}{j}\frac{j}{n}=\delta _{i,1}$
and $A_{r}\overset{d}{\rightarrow }1$ as $r\uparrow \infty ,$ the degenerate
state reached when the most recent common ancestor (MRCA) is attained.

$\left( ii\right) $ Non-degenerate solutions of $A_{\infty }$ are obtained
when considering bias mechanisms with $p\left( 0\right) =0$ and $p\left(
1\right) =1$.

$\left( iii\right) $ Consider any biased WF model with $p\left( 0\right) =0$
and $p\left( 1\right) =1$ for which $p\left( x\right) \underset{x\uparrow 0}{%
\sim }$ $\lambda x$, $\lambda >1.$ Then, due to large sample asymptotic
independence: 
\begin{equation*}
\mathbf{\nu }_{n}\overset{d}{\rightarrow }\mathbf{\xi }_{\infty },
\end{equation*}
where $\mathbf{\xi }_{\infty }$ is an iid sequence with $\xi _{1}\overset{d}{%
\sim }$ Poisson$\left( \lambda \right) $ (as it can easily be checked by the
Poisson limit to the binomial distribution). In this case, the limiting
extinction probability of $\left( N_{r};r\in \Bbb{N}\right) $ given $N_{0}=m$
is $\lim_{n\uparrow \infty }\rho _{\infty }\left( m\right) =\rho ^{m}$, $%
m=1,2,...$, where $0<\rho <1$ is the smallest solution to the fixed point
equation 
\begin{equation*}
x=e^{-\lambda \left( 1-x\right) }.
\end{equation*}
$\rho $ is the singleton extinction probability of a super-critical
Galton-Watson process with offspring distribution Poisson$\left( \lambda
\right) $. More precisely, proceeding as in M\"{o}hle '$1994$, Theorem $4.5$%
, we have 
\begin{equation*}
n\left( \rho ^{m}-\rho _{\infty }\left( m\right) \right) =\rho ^{m}\left( 
\frac{1-\rho }{1+\lambda \rho }m^{2}+\frac{\lambda \left( 1-\rho \right)
\rho }{1-\left( \lambda \rho \right) ^{2}}m\right) ,
\end{equation*}
showing that the convergence of $\rho _{\infty }\left( m\right) $ to $\rho
^{m}$ is of order $n^{-1}$. As a result, we get the asymptotic normality: 
\begin{equation*}
\frac{1}{\sqrt{n}}\left( A_{\infty }-n\left( 1-\rho \right) \right) \overset{%
d}{\underset{n\uparrow \infty }{\rightarrow }}\mathcal{N}\left( 0,\frac{\rho
\left( 1-\rho \right) }{1+\lambda \rho }\right) .
\end{equation*}

Intuitively, $\frac{1}{n}\Bbb{E}\left[ n-A_{\infty }\right] =\rho _{\infty
}\left( 1\right) =\Bbb{P}_{1}\left( N_{\infty }=0\right) \rightarrow \rho $,
showing that $\Bbb{E}\left[ A_{\infty }\right] \underset{n\uparrow \infty }{%
\sim }n\left( 1-\rho \right) $ and 
\begin{equation*}
\frac{1}{n\left( n-1\right) }\Bbb{E}\left[ \left( n-A_{\infty }\right)
\left( n-1-A_{\infty }\right) \right] =\rho _{\infty }\left( 2\right) =\Bbb{P%
}_{2}\left( N_{\infty }=0\right) \rightarrow \rho ^{2},
\end{equation*}
showing that $\sigma ^{2}\left( A_{\infty }\right) \underset{n\uparrow
\infty }{\sim }n\rho \left( 1-\rho \right) /\left( 1+\lambda \rho \right) .$

For the quadratic example 3.3, $p\left( x\right) =x\left( 1+c-cx\right) $,
with $c\in \left[ 0,1\right] $, $\lambda =1+c>1$ as soon as $c>0$. When $%
c\in \left( 0,1\right] $, we thus always have asymptotic normality. For the
example 3.1 with selection, $p\left( x\right) =\left( 1+s\right) x/\left(
1+sx\right) ,$ with $s>-1:$ $p\left( x\right) \underset{x\uparrow 0}{\sim }%
\left( 1+s\right) x$ and so $\lambda =1+s$. We have asymptotic normality
only when $s>0$, i.e. when the fitness is advantageous (corresponding as
required to complete monotonicity of corresponding $q=1-p$).\newline

Note that this asymptotic behavior does not hold for the Lipshitz-continuous
admissible mechanism $p\left( x\right) =x^{\gamma }$ of example 5.1.1\textbf{%
\ }(or more generally for mechanisms satisfying $p\left( x\right) \underset{%
x\downarrow 0}{\sim }\theta x^{\gamma }$, $\theta >0$ as in the compound
bias example\textbf{\ }5.3.5) with $\gamma \in \left( 0,1\right) $ because
its behavior near $0$ is not linear. This puzzling class of models seems to
deserve a special study as deviation to normality is expected. We postpone
it to a future work.

\section{Concluding remarks}

In this Note, we focused on discrete non-neutral Wright-Fisher models and on
the conditions on the bias probabilities under which the forward branching
dynamics is amenable to a dual discrete ancestral coalescent. It was shown
that it concerns a large class of non-neutral models involving completely
monotone bias probabilities. Several examples were supplied, some standard,
some less classical. The Wright-Fisher model with forward binomial
transition matrix is a particular instance of the Dirichlet model with
Dirichlet-binomial transition matrix. Following the same lines, using the
representation of the Dirichlet binomial distribution as a beta mixture of
the binomial distribution, it would be interesting to exhibit the
corresponding conditions on the bias mechanism, were the starting point to
be a forward Dirichlet branching process. Also of particular interest in
this respect would be the discrete non-neutral Moran models whose forward
transition matrices are simpler because of their tridiagonal Jacobi
structure. We hope to be able to consider shortly these cases (and maybe
others) in a future work.

\end{document}